\documentclass[12pt]{article} \usepackage{amsfonts}
\usepackage{amssymb}
\begin{document}
\baselineskip=14pt
\pagestyle{plain}
{\Large
\newcommand{\mun}{\mu_{n,j}}
\newcommand{\wnj}{W_{N,j}(\mu)}
\newcommand{\ntm}{|2n-\theta-\mu|}
\newcommand{\pwp}{PW_\pi^-}
\newcommand{\cpm}{\cos\pi\mu}
\newcommand{\spm}{\frac{\sin\pi\mu}{\mu}}
\newcommand{\tet}{(-1)^{\theta+1}}
\newcommand{\dm}{\Delta(\mu)}
\newcommand{\dl}{\Delta(\lambda)}

\newcommand{\dnm}{\Delta_N(\mu)}
\newcommand{\sni}{\sum_{n=N+1}^\infty}
\newcommand{\muk}{\sqrt{\mu^2+q_0}}
\newcommand{\agt}{\alpha, \gamma, \theta,}
\newcommand{\muq}{\sqrt{\mu^2-q_0}}
\newcommand{\smn}{\sum_{n=1}^\infty}
\newcommand{\lop}{{L_2(0,2\pi)}}
\newcommand{\tdm}{\tilde\Delta_+(\mu_n)}
\newcommand{\dmn}{\Delta_+(\mu_n)}
\newcommand{\dd}{D_+(\mu_n)}
\newcommand{\ddd}{\tilde D_+(\mu_n)}
\newcommand{\sn}{\sum_{n=1}^\infty}
\newcommand{\pnn}{\prod_{n=1}^\infty}
\newcommand{\aln}{\alpha_n(\mu)}
\newcommand{\emp}{e^{\pi|Im\mu|}}
\newcommand{\ppp}{\prod_{p=p_0}^\infty}
\newcommand{\ppk}{\sum_{p=p_0}^{\infty}\sum_{k=1}^{[\ln p]}}
\newcommand{\pps}{\sum_{p=p_0}^\infty}

 \medskip
\medskip

\centerline {\bf Characterization of the spectrum of irregular boundary value}
\centerline {\bf problem for the Sturm-Liouville operator}

\medskip
\medskip
\centerline { Alexander Makin}
\medskip
\medskip

\begin{quote} Abstract. We consider the spectral problem generated by the Sturm-Liouville equation with an arbitrary
complex-valued potential $q(x)\in L_2(0,\pi)$
and irregular boundary conditions. We establish necessary and sufficient conditions
for a set of complex numbers to be the spectrum
of such an operator.
\end{quote}

\medskip
\medskip

In the present paper, we consider the eigenvalue problem for the Sturm-Liouvulle equation
$$
u''-q(x)u+\lambda u=0 \eqno (1)
$$
on the interval
 $(0,\pi)$
with the boundary conditions
$$
u'(0)+(-1)^\theta u'(\pi)+bu(\pi)=0,\quad
u(0)+(-1)^{\theta+1}u(\pi)=0,
 \eqno (2)
$$
where $b$ is a complex number, $\theta=0,1$, and the function $q(x)$ is an arbitrary
complex-valued function of the class
$L_2(0,\pi)$.

Denote by $c(x,\mu), s(x,\mu)$ $(\lambda=\mu^2)$ the fundamental system of solutions to
(1) with the initial conditions $c(0,\mu)=s'(0,\mu)=1$, $c'(0,\mu)=s(0,\mu)=0$.
The following identity is well known
$$
c(x,\mu)s'(x,\mu)-c'(x,\mu)s(x,\mu)=1.\eqno(3)
$$
Simple calculations show that the characteristic equation of
(1), (2) can be reduced to the form $\Delta(\mu)=0$,
where
$$
\Delta(\mu)=c(\pi,\mu)-s'(\pi,\mu)+(-1)^{\theta+1}bs(\pi,\mu). \eqno(4)
$$
The characteristic determinant $\Delta(\mu)$
of problem (1), (2), given by (4), is referred to as the characteristic determinant
corresponding to the triple $(b, \theta, q(x))$. Throughout the following
the symbol $||f||$ stands for $||f||_{L_2(0,\pi)}$, $<q>=\frac{1}{\pi}\int_0^\pi q(x)dx$.
By $\Gamma(z,r)$ we denote the disk of radius $r$ centered at a point $z$.
By $PW_\sigma$ we denote the class of entire functions $f(z)$ of exponential type
$\le\sigma$ such that $||f(z)||_{L_2(R)}<\infty$, and by $PW_\sigma^-$
we denote the set of odd functions in $PW_\sigma$.

The following two assertions provide necessary and sufficient conditions to be
satisfied by the characteristic determinant $\Delta(\mu)$.

{\bf Theorem 1.} {\it If a function $\Delta(\mu)$ is the characteristic determinant
corresponding to the triple $(b, \theta, q(x))$, then
$$
\Delta(\mu)=(-1)^{\theta+1}b\frac{\sin\pi\mu}{\mu}+\frac{f(\mu)}{\mu},
$$
where $f(\mu)\in PW_\pi^-$.}

{\bf Proof.} Let $e(x,\mu)$ be a solution to (1) satisfying the initial conditions
$e(0,\mu)=1$, $e'(0,\mu)=i\mu$, and let $K(x,t)$, $K^+(x,t)=K(x,t)+K(x,-t)$, and
$K^-(x,t)=K(x,t)-K(x,-t)$ be the transformation kernels [1] that realize the representations
$$
e(x,\mu)=e^{i\mu x}+\int_{-x}^xK(x,t)e^{i\mu t}dt,
$$
$$
c(x,\mu)=\cos\mu x+\int_0^xK^+(x,t)\cos\mu tdt,
$$
$$
s(x,\mu)=\frac{\sin\mu x}{\mu}+\int_0^xK^-(x,t)\frac{\sin\mu t}{\mu}dt.\eqno(5)
$$
It was shown in [2] that
$$
c(\pi,\mu)=\cos\pi\mu+\frac{\pi}{2}<q>\spm-\int_0^\pi\frac{\partial K^+(\pi,t)}{\partial t}\frac{\sin\mu t}{\mu}dt,\eqno(6)
$$

$$
s'(\pi,\mu)=\cos\pi\mu+\frac{\pi}{2}<q>\spm+\int_0^\pi\frac{\partial K^-(\pi,t)}{\partial x}\frac{\sin\mu t}{\mu}dt.\eqno(7)
$$
Substituting the right-hand sides of expressions (5), (6), (7) into (4), we obtain
$$
\begin{array}{c}
\dm=(-1)^{\theta+1}b\spm+\\
+\frac{1}{\mu}\int_0^\pi[-\frac{\partial K^+(\pi,t)}{\partial t}-\frac{\partial K^-(\pi,t)}{\partial x}+(-1)^{\theta+1}bK^-(\pi,t)]\sin\mu tdt.
\end{array}
$$
This relation, together with the Paley-Wiener theorem implies the assertion of
Theorem 1.

{\bf Theorem 2.} {\it Let a function $u(\mu)$ have the form
$$
u(\mu)=(-1)^{\theta+1}b\frac{\sin\pi\mu}{\mu}+\frac{f(\mu)}{\mu},\eqno(8)
$$
where $f(\mu)\in PW_\pi^-$, $b$ is a complex number.
Then, there exists a function $q(x)\in L_2(0,\pi)$ such that the characteristic
determinant correponding to the triple $(b, \theta, q(x))$ satisfies $\Delta(\mu)=u(\mu)$.}

{\bf Proof.} Since [3]
$$|f(\mu)|\le C_1||f(\mu)||_{L_2(R)}e^{\pi|Im\mu|},\eqno(9)$$
it follows that there exists an arbitrary large positive integer
$N$ such that
$$
|u(\mu)|<1/10\eqno(10)
$$
on the set $|Im\mu|\le1$, $Re\mu\ge N$.
Let $\mu_n$ $(n=1,2,\ldots)$ be a strictly monotone increasing sequence of positive
numbers such that $|\mu_n-(N+1/2)|<1/10$ if $1\le n\le N$
and $\mu_n=n$ if $n\ge N+1$. Consider the function
$$
s(\mu)=\pi\prod_{n=1}^\infty\frac{\mu_n^2-\mu^2}{n^2}=\spm\prod_{n=1}^N\frac{\mu_n^2-\mu^2}{n^2-\mu^2}.\eqno(11)
$$
Obviously, all zeros of the function $s(\mu)$ are simple, and, in addition, the
inequality
$$
(-1)^n\dot s(\mu_n)>0.\eqno(12)
$$
holds for any
$n$.
It was shown in [4] that
$$
\dot s(n)=\frac{\pi(-1)^n}{n}(1+C_0n^{-2}+O(n^{-4})),\eqno(13)
$$
where $C_0$ is some constant, and the asymptotic formula
$$
s(\mu)=\spm+O(\mu^{-3})\eqno(14)
$$
holds in the strip $|Im\mu|\le1$.

Consider the equation
$$
z^2-u(\mu_n)z-1=0.\eqno(15)
$$
It has the roots
$$
c_n^\pm=\frac{u(\mu_n)\pm\sqrt{u^2(\mu_n)+4}}{2}.\eqno(16)
$$
It follows from (10) that for any $n$ all numbers $c_n^+$ lie in the disk $\Gamma(1,1/2)$ and
all numbers
$c_n^-$ lie in the disk $\Gamma(-1,1/2)$. Let for even
$n$ $c_n=c_n^+$, and for odd $n$ $c_n=c_n^-$.
Then $(-1)^nRec_n>0$
 for any $n=1, 2,\ldots$.
This, together with (12) implies that
$Rew_n>0$ for any $n$, where
$$
w_n=\frac{c_n}{\mu_n\dot s(\mu_n)}.\eqno(17)
$$

We set $F(x,t)=F_0(x,t)+\hat F(x,t)$, where
$$
F_0(x,t)=\sum_{n=1}^N\left(\frac{2c_n}{\mu_n\dot s(\mu_n)}
\sin\mu_nx\sin\mu_nt-
\frac{2}{\pi}\sin nx\sin nt\right),
$$
$$
\hat F(x,t)=\sum_{n=N+1}^\infty\left(\frac{2c_n}{\mu_n\dot s(\mu_n)}
\sin\mu_nx\sin\mu_nt-
\frac{2}{\pi}\sin nx\sin nt\right).\eqno(18)
$$
One can readily see that $F_0(x,t)\in C^\infty(R^2)$.
Consider the function $\hat F(x,t)$.
If $n\ge N+1$, then, by taking into account
(9), (16) and the rule for choosing the roots of equation (15), we obtain
$$
c_n=(-1)^n+\frac{f(n)}{2n}+O(1/n^2).
\eqno(19)
$$
It follows from (9), (13), (18), and (19) that
$$
\begin{array}{c}
\hat F(x,t)=\sum_{n=N+1}^\infty\frac{2}{\pi}\left(\frac{1+(-1)^n\frac{f(n)}{2n}+O(1/n^2)}{1+c_0/n^2+O(1/n^4)}-1
\right)\sin nx\sin nt=\\
=\sum_{n=N+1}^\infty\frac{2}{\pi}[(1+(-1)^n\frac{f(n)}{2n}+O(1/n^2))\times\\(1-c_0/n^2+O(1/n^4))-1
]\sin nx\sin nt=\\
=\frac{2}{\pi}\sum_{n=N+1}^\infty((-1)^{n}\frac{f(n)}{2n}+O(1/n^2))\sin nx\sin nt=\\
=(\hat G(x-t)-\hat G(x+t))/2,

\end{array}
$$
where $$\hat G(y)=\frac{2}{\pi}\sum_{n=N+1}^\infty((-1)^{n}\frac{f(n)}{2n}+O(1/n^2))\cos ny.$$
The relation
$$
\sum_{n=1}^\infty|f(n)|^2=\frac{1}{2}||f(\mu)||_{L_2(R)},
$$ which follows from the Paley-Wiener theorem, together with the Parseval equality, implies
that $\hat G(y)\in W_2^1[0,2\pi]$. Therefore, we obtain the representation
$$
F(x,t)=F_0(x,t)+(\hat G(x-t)-\hat G(x+t))/2,\eqno(20)
$$
where the functions $F_0(x,t)$ and $\hat G(y)$ belong to the above-mentioned classes.

Now let us consider the Gelfand-Levitan equation
$$
K(x,t)+F(x,t)+\int_0^xK(x,s)F(s,t)ds=0\eqno(21)
$$
and prove that it has a unique solution in the space $L_2(0,x)$ for each $x\in[0,\pi]$.
To this end, it suffices  to show that the corresponding homogeneous equation
has only the trivial solution.

Let $f(t)\in L_2(0,x)$. Consider the equation
$$
f(t)+\int_0^xF(s,t)f(s)ds=0.
$$
Following [4], by multiplying the last equation by $\bar f(t)$ and by integrating the
resulting relation over the interval
$[0,x]$, we obtain
$$\begin{array}{c}
\int_0^x|f(t)|^2dt+
\sum_{n=1}^\infty\frac{2c_n}{\mu_n\dot s(\mu_n)}\int_0^x\bar f(t)\sin\mu_ntdt\int_0^xf(s)\sin\mu_nsds-\\
-\sum_{n=1}^\infty\frac{2}{\pi}\int_0^x\bar f(t)\sin ntdt\int_0^xf(s)\sin nsds=0.
\end{array}
$$
This, together with the Parseval equality for the function system $\{\sin nt\}_1^\infty$
on the interval $[0,\pi]$ implies that
$$
\sum_{n=1}^\infty w_n|\int_0^xf(t)\sin\mu_ntdt|^2=0,
$$
where the $w_n$ are the numbers given by (17). Since $Rew_n>0$, we see that
$\int_0^xf(t)\sin\mu_ntdt=0$ for any $n=1,2,\ldots$. Since [5, 6]
the system $\{\sin\mu_nt\}_1^\infty$ is complete on the interval $[0,\pi]$, we have $f(t)\equiv0$ on
$[0,x]$.

Let $\hat K(x,t)$ be a solution of equation (21), and let $\hat q(x)=2\frac{d}{dx}\hat K(x,x)$;
then it follows [4] from (20) that $\hat q(x)\in L_2(0,\pi)$. By $\hat s(x,\mu)$, $\hat c(x,\mu)$
we denote the fundamental solution system of equation (1) with potential $\hat q(x)$ and
the initial conditions $\hat s(0,\mu)=\hat c'(0,\mu)=0$, $\hat c(0,\mu)=\hat s'(0,\mu)=1$.
By reproducing the corresponding considerations in [4], we obtain $\hat s(\pi,\mu)\equiv s(\mu)$, whence
it follows that the numbers $\mu_n^2$ form the spectrum of the Dirichlet problem for equation (1) with
potential $\hat q(x)$, and $\hat c(\pi,\mu_n)=c_n$, which, together with
identity (3), implies that
$\hat s'(\pi,\mu_n)=1/c_n$.

Let $\hat\Delta(\mu)$ be the characteristic determinant
corresponding to the triple $(b,\theta,\hat q(x))$. Let us prove
that $\hat\Delta(\mu)\equiv u(\mu)$. By Theorem 1, the function
$\hat \Delta(\mu)$ admits the representation
$$
\hat\Delta(\mu)=
(-1)^{\theta+1}b\frac{\sin\pi\mu}{\mu}+\frac{\hat f(\mu)}{\mu},
$$
where $\hat f(\mu)\in PW_\pi^-$.
By taking into account (4) and the fact that the numbers $c_n$ are roots of equation (15),
we have
$$
\hat\Delta(\mu_n)=
\hat c(\pi,\mu_n)-\hat s'(\pi,\mu_n)+(-1)^{\theta+1}b\hat s(\pi,\mu_n)
=c_n-c_n^{-1}=u(\mu_n).
$$
Hence it follows that the function
$$
\Phi(\mu)=\frac{u(\mu)-\hat\Delta(\mu)}{s(\mu)}=\frac{f(\mu)-\hat f(\mu)}{\mu s(\mu)}
$$
is an entire function on the complex plane. Since the function $g(\mu)=f(\mu)-\hat f(\mu)$
belongs to $\pwp$,
it follows from (9) that
$$
|g(\mu)|\le C_2e^{\pi|Im\mu|}.\eqno(22)
$$
Relation (11) implies that if $|Im\mu|\ge1$, then $$
|\mu s(\mu)|\ge C_3e^{\pi|Im\mu|}\eqno(23)
$$
$(C_3>0)$. hence we obtain $|Im\mu|\ge1$ if $|\Phi(\mu)|\le C_2/C_3$.

By $H$ we denote the union of the vertical segments $\{z: |Rez|=n+1/2, |Imz|\le1\}$,
where $n=N+1,N+2,\ldots$. It follows from (11) that if $\mu\in H$, then
$|\mu s(\mu)|\ge C_4>0$. The last inequality, together with (22), (23), and the
maximum principle for the absolute value of an analytic function, implies that
 $|\Phi(\mu)|\le C_5$ in the strip $|Im\mu|\le1$.
Consequently, the function $\Phi(\mu)$ is bounded on the entire complex plane
and hence identically constant by the Liouville theorem. It follows from the Paley-Wiener theorem and the
Riemann lemma [1] that if $|Im\mu|=1$, then $\lim_{|\mu|\to\infty}g(\mu)=0$,
whence we obtain $\Phi(\mu)\equiv0$.

The proof of Theorem 2 is complete.

Further we consider problem (1), (2) under the supplementary condition $b\ne0$.

{\bf Theorem 3.} {\it For a set $\Lambda$ of complex numbers to be the spectrum of problem (1), (2) it is
necessary and sufficient that it has the form $\Lambda=\{\lambda_n\}$, where $\lambda_n=\mu_n^2$,
$$
\mu_n=n+r_n,
$$
where
 $\{r_n\}\in l_2$,
 $n=1,2,\ldots$.
}

Necessity. It follows from Theorem 1 that the characteristic
equation of problem (1), (2) can be reduced to the form
$$
(-1)^\theta b\frac{\sin\pi\mu}{\mu}=\frac{f(\mu)}{\mu},\eqno(24)
$$
where $f(\mu)\in\pwp$. It was shown in [1] that equation (24) has the roots
$\mu_n=n+r_n$,
 where $r_n=o(1)$,
 $n=1,2,\ldots$.
 Hence it follows that
$$
\sin\pi r_n=(-1)^{\theta+n}f(n+r_n)/b.
$$
Since $\{f(n+r_n)\}\in l_2$, by [4], it follows that $\{r_n\}\in l_2$.

Sufficiency. Let the set $\Lambda$ admit the representation of the above-mentioned form.
We denote
$$
u(\mu)=\tet b\pi\mu\pnn\left(\frac{\lambda_n-\mu^2}{n^2}\right).
$$
It follows from [10] and the conditions of the theorem that the infinite product in the right-hand side
of the last equality converges uniformly in any bounded domain. We denote
$\phi(\mu)=\tet b\sin\pi\mu-u(\mu)$.
Let us prove that $\phi(\mu)\in\pwp$. Evidently, $\phi(\mu)$ is an odd entire
function. Obviously,
$$
|Im\mu_n|\le M,
$$
where $M$ is some constant. By $\Gamma$ we denote the union of disks of radius
$1/2$ centered at the points $n=1,2,\ldots$. Let $\mu\notin\Gamma$;
then it follows from well known identity
$$
\sin\pi\mu=\pi\mu\pnn\frac{n^2-\mu^2}{n^2}\eqno(25)
$$
that
$$
\phi(\mu)=\sin\pi\mu(1-\phi_0(\mu)),\eqno(26)
$$
where
$$
\phi_0(\mu)=\pnn(1+\alpha_n(\mu)),\eqno(27)
$$
where
$$
\alpha_n(\mu)=\frac{\mu_n^2-n^2}{n^2-\mu^2}.\eqno(28)
$$
Let us study the function $\phi_0(\mu)$ for $Re\mu\ge0$, $\mu\notin\Gamma$. One can
readily see that
$$
\begin{array}{c}
|\mu_n+n|<cn,\quad |n+\mu|>|\mu|, \\
 |n-\mu|>1/2,
\quad |n+\mu|>n.
\end{array}
\eqno(29)
$$

It follows from (28) and (29) that
$$
\begin{array}{c}
\sn|\aln|\le c_1\sn\frac{|r_n|n}{|n+\mu||n-\mu|}\le\\
\le c_1\sn\frac{|r_n|}{|n-\mu|}\le c_1(\sn\frac{|r_n|^2}{|n-\mu|^{1/2}}+\sn\frac{1}{|n-\mu|^{3/2}})\le c_2.
\end{array}\eqno(30)
$$
Obviously, the inequality
$|r_n|<1/(8c_1)$ holds
for all $n>N$, where $N$ is a sufficiently large number.
 It follows from (28) and (29) that for all $n>N$
$$
|\aln|\le1/4.\eqno(31)
$$

Relations (30), (31) and the obvious inequality
$$
|\ln(1+z)|\le2|z|,\eqno(32)
$$
valid for $|z|\le1/4$ imply that
$$
\sum_{n=N+1}^\infty|\ln(1+\aln)|\le c_3,
$$
moreover, here and throughout the following, we choose the branch of $\ln(1+z)$ that is zero
for $z=0$.
Now, by [7],
$$
\prod_{n=N+1}^\infty|1+\aln|\le e^{c_3},
$$
consequently,
$$
\prod_{n=1}^\infty|1+\aln|\le c_4e^{c_3}.\eqno(33)
$$
Since $\phi_(\mu)$ is an even function, it follows from (26), (27), and (33) that
$$
|\phi(\mu)|\le c_5\emp.\eqno(34)
$$
The maximum principle implies that inequality (34)
is valid in the entire complex plane, therefore, the function $\phi(\mu)$ is an entire function
function of exponential type $\le\pi$.
Evidently, $|r_n|<c_6$. Notice, that if $|Im\mu|\ge\tilde C$, where $\tilde C=4c_1c_6$,
then estimate (31) holds for any $n=1,2,\ldots$.
In the domain $|Im\mu|\ge\tilde C$
we define a function
$$
W(\mu)=\ln\phi_0(\mu)=
\sn\ln(1+\aln),
$$
then we have
$$
\phi(\mu)=\sin\pi\mu(1-e^{W(\mu)}).\eqno(35)
$$
Let us estimate the function $W(\mu)$. One can readily see that
$$
\lim_{|Im\mu|\to\infty}\left(\sn\frac{|r_n|^2}{|n-\mu|^{1/2}}+\sn\frac{1}{|n-\mu|^{3/2}}\right)=0.
$$
This, together with (30) and (32) implies that
the inequality
$$
\begin{array}{c}
|W(\mu)|\le
\sn|\ln(1+\aln)|\le 2
\sn|\aln|)\le1/4
\end{array}
$$
holds in the domain $|Im\mu|\ge \hat C$, where $\hat C$ is a sufficiently large
number.
Therefore, it follows from the elementary inequality $|1-e^z|\le2|z|$, valid for $|z|\le1/4$, that
$|1-e^{W(\mu)}|\le2|W(\mu)|$, which, together with (35) implies that
$$
|\phi(\mu)|\le c_7|W(\mu)|.\eqno(36)
$$
for $\mu\in l$,
where $l$ is the line $Im\mu=\hat C$.
Let us prove the inequality
$$
\int_l|W(\mu)|^2d\mu<\infty.\eqno(37)
$$
Let $\mu\in l^+$, where $l^+$ is the ray $Im\mu=\hat C,
Re\mu\ge0$. The elementary inequality
$$
|\ln(1+z)-z|\le|z|^2,
$$
valid for $|z|\le1/2$, implies that

$$
|W(\mu)|\le |S_1(\mu)|+S_2(\mu),
$$ where
$$
S_1(\mu)=\sn\aln, \quad S_2(\mu)=\sn|\aln|^2.
$$
Set
$$
I_m=\int_{l^+}|S_m(\mu)|^2d\mu
$$
$(m=1,2)$. First, consider the integral $I_1$. One can readily see that
$$
\begin{array}{c}
I_1=\int_{l^+}\left|\sn\frac{(\mu_n+n)r_n}{(n+\mu)(n-\mu)}\right|^2d\mu
\le c_7(
\int_{l^+}\left|\mu\sn\frac{r_n}{(n+\mu)(n-\mu)}\right|^2d\mu+\\
+\int_{l^+}\left|\sn\frac{r_n}{(n+\mu)}\right|^2d\mu
+\int_{l^+}\left|\sn\frac{r_n^2}{(n+\mu)(n-\mu)}\right|^2d\mu).
\end{array}\eqno(38)
$$
The convergence of the first and the third integrals in the right-hand side of (38) was establised in [8],
By [9] so is the second integral.
It is readily seen that
$$
I_2\le c_8
\int_{l^+}\left|\sn\frac{|r_n|^2}{|n-\mu|^2}\right|^2d\mu.\eqno(39)
$$
By [8] the integral in right-hand side of (39) is convergent.

From the last inequality, the convergence of the integral $I_1$ and the evenness of the function $W(\mu)$,
we find that inequality (37) is valid. It follows from (36), (37) and
[3] that
$$
\int_R|\phi(\mu)|^2d\mu<\infty,
$$
consequently, $\phi(\mu)\in\pwp$, which, together with Theorem 2,
proves theorem 3.

Consider problem (1), (2) if $b=0$. Substituting the functions $c(x,\mu)$, $s(x,\mu)$
into boundary conditions and taking into account (3), we find that each root
subspace contains one eigenfunction and, possibly, associated functions. The characteristic
equation has the form
$$
\frac{f(\mu)}{\mu}=0,
$$
where $f(\mu)\in PW_\pi^-$. Let us consider two examples.

1). Set
$$
f_1(\mu)=\frac{\sin^k(\alpha\pi\mu/k)\sin^k((1-\alpha)\pi\mu/k)}{\mu^{2k-1}},
$$
where $k$ is an arbitrary natural number, and $\alpha$ is an
irrational number, $0<\alpha<1$. Obviously, $f_1(\mu)\in
PW_\pi^-$. Then, by Theorem 2, there exists a potential $q_1(x)\in
L_2(0,\pi)$, such that the corresponding characteristic
determinant $\Delta_1(\mu)=f_1(\mu)/\mu$. Since the equations
$\sin(\alpha\pi\mu/k)=0$ and $\sin((1-\alpha)\pi\mu/k)=0$ have no
common roots, except zero, we see that each root subspace of
problem (1), (2) with potential $q_1(x)$ contains one
eigenfunction and associated functions up to order $k-1$. One can
readily see that $|\Delta_1(\mu)|\ge ce^{|Im\mu|\pi}|\mu|^{1-2k}$
$(c>0)$, if $\mu$ belongs to a sequence of infinitely expanding
contours. Then, by [11], the system of eigen-and associated
functions of problem (1), (2) is complete in $L_2(0,\pi)$.

2). Set $f(\mu)=\sin^2(\pi\mu/2)/\mu$. It follows from (25) that
$$
f(\mu)=\frac{\pi^2}{4}\mu\pnn\left(\frac{(2n)^2-\mu^2}{(2n)^2}\right)^2.
$$
We denote
$$
u(\mu)=\frac{\pi^2}{4}\mu\pnn\left(\frac{\mu_n^2-\mu^2}{(2n)^2}\right)^2,\eqno(40)
$$
where
$\mu_n=2n$, if $n\ne2^p+k$, $k=1,\ldots,[\ln p]$, $p=p_0,p_0+1,\ldots$.
$\mu_n=2^{p+1}$, if $n=2^p+k$, $k=1,\ldots,[\ln p]$, $p=p_0,p_0+1,\ldots$ $(p_0\ge10)$.
It can easily be checked that
$$
\lim_{n\to\infty}\frac{\mu_n}{n}=2,\quad 0<c_1<\pnn\frac{\mu_n}{n}<\infty.
$$
This, together with [7] implies that the infinite product in right-hand side of
(40) uniformly convergents in any bounded domain of the complex plane,
therefore, $u(\mu)$ is an entire analytical function.

Let us prove that
$u(\mu)\in PW_\pi^-$. We set $\psi(\mu)=f(\mu)-u(\mu)$.
By $\Gamma$ we denote the union of disks of radius $1$ centered at the points
$2n$, $n=1,2,\ldots$. Let $\mu\notin\Gamma$,
$Re\mu\ge0$; then
$$
\psi(\mu)=f(\mu)(1-\phi(\mu)),
$$
where
$$
\begin{array}{c}
\phi(\mu)=(\ppp A_p(\mu))^2,\quad A_p(\mu)=\prod_{k=1}^{[\ln p]}(1+\alpha_{p,k}(\mu)),
\\ \alpha_{p,k}(\mu)=\frac{-(2^{p+1}+k)k}{(2^p+k-\mu/2)(2^p+k+\mu/2)}.
\end{array}
$$
Trivially,
$$
|\alpha_{p,k}(\mu)|\le2\ln p.\eqno(41)
$$

Consider two cases. Let $|\mu/2-2^p|\ge2^p/10$ for all $p=p_0, p_0+1,\ldots$;
then
$$
|\alpha_{p,k}(\mu)|\le\frac{2\ln p}{|2^p+k-\mu/2|}\le1/4.\eqno(42)
$$
Consider the function
$$
F(\mu)=\ppk \ln(1+\alpha_{p,k}(\mu)).
$$
It follows from (31), (41), and (42)
that
$$
\begin{array}{c}
|F(\mu)|=\ppk|\ln(1+\alpha_{p,k}(\mu))|
\le\ppk|\alpha_{p,k}(\mu))|\le\\
\le c_1\pps\frac{\ln^2 p}{|2^p-\mu/2|}\le c_2\pps\frac{\ln^2 p}{2^p}\le c_3.
\end{array}
$$
Now, by [7],
$$
|\phi(\mu)|\le e^{2c_3}.\eqno(43)
$$

Suppose for some $\tilde p\ge p_0$
$$
|\mu/2-2^{\tilde p}|<2^{\tilde p}/10.\eqno(44)
$$
Evidently, $\phi(\mu)=\gamma_{\tilde p}^2 (\mu)\beta_{\tilde p}^2 (\mu)$,
where
$$
\gamma_{\tilde p}(\mu)=\prod_{k=1}^{[\ln p]}(1+\alpha_{\tilde p,k}(\mu)),\quad
\beta_{\tilde p}(\mu)=\prod_{p=p_0, p\ne\tilde p}^{\infty}\prod_{k=1}^p (1+\alpha_{p,k}(\mu))
$$
Arguing as above, we see that
$$
|\beta_{\tilde p}(\mu)|\le c_4,
$$
where $c_4$ does not depend of $\tilde p$. It follows from (41) and (44) that
$$
|\gamma_{\tilde p}|^2\le(1+2\ln\tilde p)^{2\ln\tilde p}\le c_52^{\tilde p/3}\le c_6|\mu|^{1/3}.
$$
From this inequality and the evenness of the function $\phi(\mu)$, we find that
the inequality
$$
|\psi(\mu)|\le |f(\mu)|(|\mu|+3)^{1/3}
$$
is valid outside
$\Gamma$.
Since
 $|f(\mu)|\le c_7/(|\mu|+3)$ in the strip
 $\Pi: |Im\mu|\le2$, we see that the relation
$$
|\psi(\mu)|\le c_8(|\mu|+3)^{-2/3}
$$
holds on the
set
$\Pi\setminus\Gamma$.
The last inequality and the maximum principle imply that
$$
|\psi(\mu)|\le c_8(|\mu|+1)^{-2/3}
$$
in the strip
 $\Pi$, consequently,
$\psi(\mu)\in PW_\pi^-$, hence,
$u(\mu)\in PW_\pi^-$.
Then by theorem 2 there exists a potential $q(x)\in L_2(0,\pi)$,
such that the correponding characteristic determinant $\dm=u(\mu)/\mu$.
This yields that the dimensions of root subspaces of problem (1), (2) with potential
$q(x)$ increase  infinitely, and the system of root functions contains associated
functions of arbitrarily high order.

By $Q_p$ we denote the disk $|\mu/2-2^p|<2^p/10$, $Q=\bigcup_{p=p_0}^\infty$,
$D=\mathbb{C}\setminus(\Gamma\bigcup Q)$. In the domain $D$ we consider the function
$\tilde\phi(\mu)=f(\mu)/u(\mu)$. It is readily seen that
$$
\begin{array}{c}
\tilde\phi(\mu)=(\ppp \tilde A_p(\mu))^2,\quad\tilde A_p(\mu)=\prod_{k=1}^{[\ln p]}(1+\tilde\alpha_{p,k}(\mu)),
\\ \tilde\alpha_{p,k}(\mu)=\frac{(2^{p+1}+k)k}{(2^p-\mu/2)(2^p+\mu/2)}.
\end{array}
$$
Arguing as above, we see that
$$
|\tilde\phi(\mu)|\le c_9.
$$
This implies that $\dm\ge c_{10}e^{|Im\mu|\pi}/|\mu|^2$ $(c_{10}>0)$
if $\mu$ belongs to a sequence of infinitely expanding contours.
 Then, by [11],
the system of eigen-and associated functions of problem (1), (2) is complete in $L_2(0,\pi)$.

Characterization of the spectrum of nonselfadjoint problem (1), (2) in the case
of separated boundary conditions was given in [12], and for periodic and antiperiodic
boundary conditions an analogous question was solved in [9].

\medskip
\medskip
\medskip

\centerline {\bf References}
\medskip
\medskip
\medskip

1. V.A. Marchenko, Sturm-Liouville Operators and Their Applications (Kiev, 1977).

2. L.A. Pastur, V.A. Tkachenko, Spectral theory of a class of one-dimensional
Schrodinger operators with limit-peroodic potentials, Trudy Moscow Mat. Obshch.
{\bf 51}, 114-168 (1988).

3. S.M. Nikolskii, Approximation of Functions of Several Variables and Embedding Theorems
(Moscow, 1977).

4. V.A. Tkachenko, Discriminants and generic spectra of nonselfadjoint Hill's operators,
 Adv. Sov. Mathem. {\bf 19}, 41-71 (1994).

5. A.M. Sedletskii, The stability of the completeness and of the minimality in $L_2$ of a
system of exponential functions, Mat. Zametki {\bf 15}, 213-219 (1974).

6. A. M. Sedletskii, Convergence of angarmonic Fourier series in systems of exponentials,
cosines and sines, Dokl. Akad. Nauk SSSR {\bf 301}, 1053-1056 (1988).

7. M. A. Lavrent'ev and B.V. Shabat, Methods of the Theory of a Complex Variable,
(Moscow, 1973).

8. A.S. Makin, Characterization of the Spectrum of Regular Boundary Value Problems for the
Sturm-Liouville Operator, Differ. Uravn., {\bf 44}, 329-335 (2008).

9. J.-J. Sansuc, V. Tkachenko, Characterization of the Periodic and Anti-periodic Spectra of
Nonselfadjoint Hill's Operators, Oper. Theory Adv. Appl. {\bf 98}, 216-224 (1997).

10. M.A. Evgrafov, Asymptotic Estimates and Entire Functions (Moscow, 1969).

11. M.M. Malamud, On Completeness on the System of Root Vector Sturm-Liouville
Operator Subject to General Boundary Conditions, Dokl. Math. {\bf 419}, 19-22 (2008).

12. S. Albeverio, R. Hryniv and Ya. Mykytyuk, On spectra of non-self-adjoint
Sturm-Liouville operators, Selecta Math. {\bf 13}, 517-599 (2008).

\medskip
\medskip
\medskip
\medskip
\medskip
\medskip

Moscow State University of Instrument-Making and Informatics, Stromynka 20,
Moscow, 107996, Russia

E-mail address: alexmakin@yandex.ru

}
\end{document}